\documentclass[11pt]{article}
\usepackage{hyperref} 
\begin{document}
\baselineskip=15pt
\title{A $p$-adic Height Function Of Cryptanalytic Significance}
\author{H. Gopalkrishna Gadiyar and R. Padma\\
AU-KBC Research Centre\\ M. I. T. Campus of Anna University\\ Chromepet, Chennai 600044, India.\\e-mail: \{gadiyar, padma\}@au-kbc.org}
\maketitle

\begin{abstract}
It is noted that an efficient algorithm for calculating a $p$-adic height could have cryptanalytic applications. 
\end{abstract}

Elliptic curves and their generalizations are an active research topic with practical applications in cryptography \cite{kob1}, \cite{miller}, \cite{kob2}. If $E$ is an elliptic curve defined over a finite field $F_p$, where $p$ is prime, and if $P$ and $Q$ are points on the curve $E$ such that $Q=nP$, then the elliptic curve discrete logarithm problem (ECDLP) is to find $n$ given $P$ and $Q$ on $E$. The essential difficulty of the ECDLP arises from the curves being defined over finite fields and the addition law having a complicated nonlinear structure. Philosophically several papers have been written where efficient algorithms have been developed by taking advantage of the fact that finite fields sit very comfortably inside the $p$-adic fields. More precisely, objects defined over finite fields can be lifted to objects over $p$-adic fields using a simple or sophisticated versions of the Hensel's lemma \cite{gadiyar:padma:maini}, \cite{finotti}. Conversely, a $p$-adic object can be truncated to give an object in a finite field. Efficient algorithms using $p$-adic techniques include the attack on anomalous curves \cite{smart}, \cite{semaev} and  \cite{satoh:araki}, $p$-adic point counting \cite{satoh}, \cite{kedlaya}, constructing elliptic curves with the required number of points \cite{broker:stevenhagen} etc. 

In this short note, we wish to bring to the attention of the cryptographic community some recent developments in the $p$-adic height function \cite{stein}, \cite{mazur:stein:tate}, \cite{wuthrich} which could be of cryptanalytic interest. 

Several years ago Abramov and Rosenbloom (unpublished handwritten manuscript) had developed a mod $p$ version of the $\sigma $- function  and used it to define a mod $p$ height. See Introduction to the paper \cite{mazur:tate} and the references therein. In \cite{abramov:rosen}, they have defined a $p$-adic $\sigma $-function $\sigma _p$ and a $p$-adic height function $h_p$. The significance of this height function is that it is quadratic (that is, $h_p(nP)=n^2h_p(P)$), unlike the local height which does not have this nice property. If a point on an elliptic curve has coordinates in the form $\left (\frac{a}{d^2}, \frac{b}{d^3}\right )$, then $h_p$ is defined as $- \log_p\left ( \frac{\sigma _p^2 }{d^2} \right )$. The transformation properties of $\sigma _p$ and $d$ are used to kill inconvenient terms in the local height so that $h_p$ becomes quadratic. As a comment on the side, we wish to point out that in Lenstra's elliptic curve factorization algorithm \cite{lenstra}, \cite{silverman:tate},  the denominator plays a key role. 

Koblitz in \cite{kob3} states that the height function protects the elliptic curve cryptosystem, the reason being essentially that rational points with very large number of digits would be required in any attack. This is true if one lifts to $Q$ but as is now folklore, the natural structure to lift to is not $Q$, but the $p$-adic field $Q_p$, from finite fields. For example, in Smart's attack on anomalous curves \cite{smart}, it is sufficient to take a $p$-fold sum and calculate modulo $p^2$. The same is true of various $p$-adic algorithms. 

Recently an efficient algorithm for computing the $p$-adic height has been invented by Mazur, Stein and Tate \cite{mazur:stein:tate}. This could possibly be used to develop an efficient algorithm to break cryptosystem based on elliptic curves.

Finally, we wish to make some historical comments. Elliptic curves over rationals $E(Q)$ have been traditionally studied because of interest in Diophantine equations and local heights studied as global heights were broken into local ones. In cryptography, as we begin with $E(F_p)$, the natural object to study becomes $E(Q_p)$. Note that the $p$-adic height $h_p$ is a mapping from $E(Q) \cap E_{1,p} \rightarrow pZ_p$, where $E_{1,p}$ consists of points in $E(Q_p)$ with Lutz filtration $\ge 1$ and this was developed in the context of the mod $p$ variant of the Birch-Swinnerton-Dyer conjecture \cite{abramov:rosen}.


\begin{thebibliography}{18}
%
\bibitem{kob1} N. Koblitz, {\it ``Elliptic Curve Cryptosystems"}, Math. Comp., {\bf 48}(1987) 203-209.
%
\bibitem{miller} V. S. Miller, {\it ``Use of elliptic curves in cryptography"}, Advances in Cryptology-CRYPTO '85(LNCS 218), 417-426, 1986. 
%
\bibitem{kob2} N. Koblitz, {\it ``Hyperelliptic Cryptosystems"}, J. Crypto., {\bf 1}(1989), 139-150.
%
\bibitem{gadiyar:padma:maini} H. Gopalkrishna Gadiyar, K M Sangeeta Maini and R. Padma, {\it ``Cryptography, Connections, Cocycles and Crystals: A $p$-adic Exploration of the Discrete Logarithm Problem"}, Progress in Cryptology - Indocrypt 2004(LNCS 3348), 305-314.
%
\bibitem{finotti} L. R. A. Finotti,{\it Canonical and Minimal Degree Liftings of Curves}, Ph. D. thesis, University of Texas at Austin, August 2001.
%
\bibitem{smart} N. P. Smart, {\it ``The discrete logarithm problem on elliptic curves of trace one"}, J. Crypto., {\bf 12}(1999), 193--196.
 
%
\bibitem{semaev} I. A. Semaev, {\it ``Evaluation of discrete logarithms on some elliptic curves"}, Math. Comp., {\bf 67}(1998), 353-356.
%
\bibitem{satoh:araki} T. Satoh and K. Araki, {\it ``Fermat quotients and the polynomial time discrete log algorithm for anomalous elliptic curves"}, Comm. Math. Univ. Sancti Pauli, {\bf 47}(1998), 81-92.
%
\bibitem{satoh} T. Satoh, {\it ``The canonical lift of an ordinary elliptic curve over a finite field and its point counting"}, J. Ramanujan
Math. Soc., {\bf 15}(2000), 247-270.
%
\bibitem{kedlaya} Kiran S. Kedlaya, {\it Counting points on hyperelliptic curves using Monsky-Washnitzer cohomology}, J. Ramanujan Math. Soc., {\bf 16}(2001), 323-338. 
%
\bibitem{broker:stevenhagen} R. Br\"{o}ker and P. Stevehagen, {\it ``Elliptic Curves with a given number of points"}, ANTS 2004, LNCS(3076), 117-131.
%
\bibitem{stein} W. Stein, {\it ``An algorithm for computing $p$-adic heights using Monsky-Washnitzer cohomology"}, Notes for a Talk at MIT on 2004-10-15, http://modular.ucsd.edu/talks/padic$\_$heights/pheight.pdf 
%
\bibitem{mazur:stein:tate} B. Mazur, W. Stein and J. Tate, {\it ``Computation of $p$-adic heights and log convergence"}, http://modular.ucsd.edu/papers/pheight
%
\bibitem{wuthrich}, C. Wuthrich, {\it The Hare and the Wolf: An essay on $p$-adic heights on elliptic curves}, Ph. D. Thesis, University of Cambridge.
%
\bibitem{mazur:tate} B. Mazur and J. Tate, {\it ``The $p$-adic sigma function"}, Duke Math. Journal {\bf 62} No.3, (1991),663-688.
%
\bibitem{abramov:rosen} S. A. Abramov and M. J. Rosenbloom, {\it ``On the Birch-Swinnerton-Dyer conjecture mod $p$"}, J. Number Theory {\bf 19}(1984), 297-300.
%
\bibitem{lenstra} H. W. Lenstra, {\it ``Factoring integers with elliptic curves"}, Annals of Math. {\bf 126}(1987), 649-673.
%
\bibitem{silverman:tate} J. H. Silverman and J. Tate, {\it ``Rational Points on Elliptic Curves"}, UTM, Springer.
%
\bibitem{kob3} N. Koblitz, {\it ``Miracles of the height function - a golden shield protecting ECC"}, ECC 2000, Essen, Germany.







\end{thebibliography}
\end{document}